\documentclass[a4paper,12pt,oneside]{article}

\usepackage{amsmath}
\usepackage{amsfonts}     
\usepackage{amstext}
\usepackage{amsthm}       
\usepackage{xspace}

\usepackage{graphicx}

\usepackage{latexsym}
\usepackage[english]{babel}
\input psfig.sty

\newtheorem{theorem}{Theorem}
\newtheorem{lemma}[theorem]{Lemma}
\newtheorem{proposition}[theorem]{Proposition}

\newtheorem{example}[theorem]{Example}
\newtheorem{definition}[theorem]{Definition}

\newtheorem{remark}[theorem]{Remark}

\begin{document}
\title{ Global and local Complexity in weakly chaotic dynamical systems\footnote{AMS: 37B99, 37B40, 37A35.} }

\author{Stefano Galatolo
\footnote{
Dipartimento di Matematica  
 Universit\`a di Pisa 
 Via Buonarroti 2/a, Pisa, Italy 
e-mail: galatolo@dm.unipi.it}}

\maketitle
\abstract {In a topological dynamical system the  complexity of an orbit is a measure of the amount of information (algorithmic information content) that is necessary to describe the orbit.
This indicator is invariant up to  topological conjugation.  
We consider this indicator of local complexity of the dynamics and
provide different examples of its behavior, showing how it can be useful
to characterize various kind of weakly chaotic dynamics.
We also provide criteria to find systems with non trivial orbit
complexity (systems where the description of the whole orbit
requires an infinite amount of information).
We consider also a global indicator of the complexity of the system.
This global indicator generalizes the topological entropy, taking
into account systems were the number of essentially different
orbits increases less than exponentially.
Then we prove that if the system is constructive (roughly speaking:
if the map can be defined up to any given accuracy using a finite amount
of information) the orbit complexity is everywhere less or equal
than the generalized topological entropy. Conversely
there are compact non constructive examples where the inequality 
is reversed, suggesting that this notion comes out  naturally
in this kind of complexity questions.}

\def\stackunder#1#2{\mathrel{\mathop{#2}\limits_{#1}}}

\tableofcontents

\section{Introduction}

 Weakly chaotic phenomena are widely studied in the physical literature.
There are connections with many physical  phenomena: self organized criticality,  the anomalous diffusion processes  and many others.

In the literature a precise definition of weak chaos is not given
and many different examples are studied. Roughly speaking a weakly chaotic
system is a system that is chaotic in some sense (for example it has sensitive
dependence to initial conditions) but it has zero entropy (KS-entropy or topological entropy) in this work we will mainly consider systems having zero topological entropy.

For the classification of weakly chaotic systems there have been proposed many 
invariants.
Some of them are defined by a generalization of  K-S entropy (\cite{blume},\cite{tsallis},\cite{takens}) orbit complexity (\cite{ionewnew}) or
based on growth rate of the number of different orbits with respect to 
a partition and an invariant measure (\cite{feren}).

Generalized orbit complexity associates to each  point an indicator
of the complexity of its orbit. This is a  measure of the growth rate of the information that is necessary to describe the orbit while time increases.
 
In this sense orbit complexity can be viewed as a local version of the entropy, where  
a global (average) notion of information (the Shannon information) is replaced by a local one 
(the algorithmic information content). 
This is also confirmed by the relation (Theorem \ref{brubno}) between entropy and orbit complexity 
that can be proved in the positive entropy case.
 
Let us consider an orbit of a discrete time dynamical system. The definition of orbit complexity associates to the orbit  
 a set of strings by a geometrical construction
and then the information content of the strings is considered to define
the complexity of the orbit. 
For this purpose we need an intrinsic (pointwise) notion of information content
of a string.

Such an intrinsic notion of information content is is given for example by the Algorithmic Information Content
( also called Kolmogorov-Chaitin complexity).
Other notions of information content can be considered (\cite{io3} \cite{CSF}) but for the purposes of this paper we will only consider AIC as a measure for the information
content that is contained in a string.

Generalized orbit complexity turns out to be related to another important
feature of chaos: the sensitivity to initial conditions.
In \cite{ionewnew} quantitative relations are proved between the growth rate of the information needed to describe an orbit, quantitative indicators of sensitivity to initial conditions and the dimension of the underlying space.

In this paper we consider a notion of orbit complexity  that is particularly
suited for systems where the information necessary to describe an orbit
increases particularly slowly as the time increases.
This give a slight modification of the orbit complexity indicators defined in \cite{ionewnew}. We also consider a sort of generalized topological entropy defining a family of invariants of topological dynamical systems that contains the classical definition of topological entropy as a special case.
This generalized topological entropy provides a family of invariants that
can distinguish between topological dynamical systems with zero
entropy, characterizing the global complexity of its behavior.

One of the main results (Theorem \ref{infinit}) of the paper is that  if a system is constructive (the map can be approximated by an algorithm, see Section 5 for a precise definition)
orbit complexity (the local indicator of complexity of the system)
is less or equal than generalized topological entropy (the global indicator)
while if the system is not constructive this inequality does not hold, proving
that constructivity comes out naturally when considering Algorithmic information content based notions of complexity.

Another main result is a criterion (Proposition \ref{mitico}) to find systems with non trivial orbit complexity. 
This criterion implies for example that a system that is chaotic in the sense
of Li and Yorke has nontrivial orbit complexity.
This criterion implies (Remark \ref{misi}) that orbit complexity provides invariants that can distinguish between dynamical systems that 
that are isomorphic in the measure preserving framework.

In Section 2 we give a short introduction to the concept of 
algorithmic information content.

In Section 3 we introduce two different notions of  complexity of single orbits
given by two different variants of the notion of information content
of a string. After the definition of this invariants of the dynamics
we give some example of its calculation in different examples of dynamical systems.
We also state an easy criterion (Theorem \ref{mitico}) to find systems with
{\em nontrivial} orbit complexity.

In Section 4 we define a {\em global} indicator of the complexity of a dynamical system. The indicator generalizes in some sense the topological entropy,
taking into account different possible asymptotic behaviors of the number 
of substantially different orbits that appears in $n$ steps of the dynamics.

In Section 5 we introduce the concept of constructivity.
Roughly speaking a map is constructive if the map can be 
approximated at any accuracy by some algorithm. A rigorous definition can
be given for maps between very general metric spaces.
Constructivity of the map underlying the dynamics is  an assumption 
that implies interesting features of orbit complexity and relation
with other indicators of complexity and chaos (section 6).

In Section 6 we prove that if the map is constructive then
 the orbit complexity of each point is less or equal than
the indicator of global complexity. Constructivity
is an essential assumption. 
An example is given to show that even in the compact case
 there are (non constructive) maps with big orbit complexity
and low global complexity.

\section{Algorithmic Information content}\label{AIC}

In this section we give a short introduction to algorithmic information theory. A more detailed exposition of algorithmic information theory can be found in \cite{Zv} or \cite{Ch}.

The AIC associates to a single string a measure of the information content of 
a string, that depends (up to a constant) only on the given string. This is a very  powerful tool and allows pointwise definitions.

Let us consider the set  $\Sigma=\{0,1\}^*$ of finite (possibly empty)  binary strings. If $s$ is a string we define $|s|$ as the length of $s$.

The Algorithmic Information Content (AIC) of a string is the length of the smallest program 
to be run by some computing machine giving the string as the output. 
In other words the AIC of a string is the length of its shorter algorithmic description. 
 
For example the algorithmic information content  
of a $2n$ bits   long periodic string  
$$s=''1010101010101010101010...''$$  is small because  
the string is output of the short program:  
 
\centerline{\em repeat $n$ times (write (``10'')).} 
 
\noindent The AIC of the string $s$  then  satisfies $AIC(s)\leq log(n)+Constant$. This is because $log (n)$  
 bits are sufficient  to  code ``$n$'' (in binary notation) and the constant represents the length of the code for the computing machine representing the instructions ``repeat...''. As it is intuitive the information content of a periodic string is very poor. 
On the other hand each $n$ bits long string $$s'=''1010110101010010110...''$$ is output of the trivial program  
$${write(''1010110101010010110...'')}.$$  
This has length $n+constant$. This implies that the AIC of each string 
is (modulo a constant which depends on the chosen computing machine) less or equal than its length.  
 
The concept of computing machine can be formalized by the theory of Turing machines or recursive functions. 
 For our intuitive approach  let us think that a computing machine is an every day computer $C$ to which it can be given some program to run. 
If we give it a program $p$ (coded by some binary string)  
to be run and the computation stops we obtain an output $s$ (another string) in this case we write 
$C(p)=s$. We can suppose that the output is a string made of digits in a finite alphabet.  $C$ then defines a  function from a subset of the set of the finite binary 
strings (where the computation stops) to the set of finite strings from a finite alphabet. In the language of theoretical  
computer science this means that $C$ defines a {\em partial recursive function}. If conversely the computation stops and the output is defined for each
input then we say that   $C$ defines a {\em total recursive function}. Recursive functions are functions whose values can be calculated by some algorithm.
By this notations we can define more formally 
 
\begin{definition}[Algorithmic Information Content] 
The Kolmogorov complexity or Algorithmic Information Content of a string  
$s$ given $C$ is the length of the smallest 
program $p$ giving $s$ as the output: 
 
$$ 
{ AIC}_C(s)=\stackunder{p\in \Sigma ,C(p)=s}{\min }|p|,  
$$ 
 
\noindent  if $s$ is not a possible output for 
the computer $C$ then ${AIC}_C(s)=\infty $ . 
\end{definition}

In the last definition the algorithmic information content of a string depends on the choice of $C$. To avoid this problem we require that $C$ is an universal machine.
Roughly speaking a computing machine is called universal if it can emulate
each other machine if an appropriate input is given.

In the examples above we have that  the programs are written  
in a "Pascal like" language and $C$ is represented by a system  
able to run such program.  $C$ is then essentially a Pascal interpret. 
 
 Let us consider $C$, the above Pascal interpret and let $D$ be a Lisp interpret.
Since using Pascal language we can write a program $L$ which is a Lisp interpret
we have that for each Lisp program $p$ we have $D(p)=C(L,p)$ and then
$|{ AIC}_C(s)-{AIC}_D (s)|\leq |L|$ for each string $s$.

A formal definition of universal computing machine of  course can be given.
In this definition we also have to specify the meaning of
the ``,'' in  ``$L,p$''. Indeed a pair of binary strings can be encoded 
into a single binary string in a way that both the strings can be recovered from the encoded string whitout loosing information. For example such an encoding can be done by adding $p$ to a self delimiting
description of $L$ (an encoding of $L$ that starts specifying its length). However we will not go into technical details, for 
our scope it is sufficient to think to universal computing
machines as our every day computer that can be programmed for general
purposes tasks.
The only important difference we have to consider is that theoretical computing
machines have virtually infinite memory, that is, while computing they can write (and then read) data
on an infinite tape.

The important property of Universal Computing Machines (UCM)
that will be used here is the following.

\begin{theorem}
If ${\cal U}$ and ${\cal U}^{\prime }$ are universal computing machines then 
\begin{equation*}
\left| {AIC}_{\cal U}(s)-{AIC}_{{\cal U}\prime}(s)\right| \leq K\left( {\cal U},{\cal U}^{\prime
}\right)
\end{equation*}
where $K\left( {\cal U},{\cal U}^{\prime }\right) $ is a constant which depends only on ${\cal U}$
and ${\cal U}^{\prime }$ but not on $s$.
\end{theorem}

This theorem sates that if we use an UCM  in the definition of  the algorithmic information content
then  this information content does not depends
on the particular machine we choose in this class up to a constant. 
Since we are interested to the asymptotic behavior of the quantity of information 
this constant is not relevant and
this remark allows to not mention the chosen machine ${\cal U}$ in the notation
${AIC}_{\cal U}(s)$ in the future.
In the remaining part of the paper the universal computing machine
that is considered in the definition of AIC will be denoted by {\cal U}.

We also  want to consider the information that is necessary to 
 reconstruct a string $s$ once another string $s'$ is known.

As it was said before there are many ways to encode a pair of strings
into a single string. 
Let us choose such an encoding $s',s\rightarrow <s',s>$ and suppose that 
it is injective and recursive.

Being universal our computing machine can be also supposed to be  able to recognize this encoding and recover both strings from the encoded string.
Now we can define 
\begin{definition} The conditional $AIC $ of $s$ given $s'$  is the 
length of the shortest program 
that is able to reconstruct the string $s$ when $s'$ is given:
$$AIC(s|s')=\stackunder {{\cal U}(<s',p>)=s}{min} |p|.$$
\end{definition}

Up to a constant the definition is independent 
on the  encoding that is chosen for the pair. This is because for each pair of chosen encodings
there is an algorithm translating one encoding in the other for each
pair of strings and this program will only add a constant in the definition
of information content.

\section{Orbit complexity}
\label{sectcom}

We give two definitions of orbit complexity. One is based on the plain  
algorithmic information content and is the definition that was given in 
\cite{ionewnew}, the other is based on the algorithmic information of a sting 
given its length. The latter is particularly suited for very regular orbits.

 Let us consider a   dynamical system $(X,T)$.  $X$ is a 
 compact metric space and $T$ is a 
  function $X\rightarrow X$. Until section 4 $T$ is not necessarily supposed to be continuous. 
Let us consider a  finite open cover  $\beta=\{B_0,B_1,...,B_{N-1}\}$ 
of $X$, that is a collection of open sets whose union is the whole space.

 We use this cover
to code the orbits of $(X,T)$  into a set of infinite strings. 
A symbolic coding of the orbits of $X$ with respect to the cover $\{B_i\} $ is a string listing the sets $B_1,..,B_n$ visited by the orbit of $x$ during 
the iterations of $T$. Since the sets $B_i$ may have non empty intersection 
then an orbit can have more than one possible coding. More precisely,
if  $x\in X$ let us define the set of symbolic orbits of $x$ with respect to $\beta$ as: 
$$\varphi _\beta(x)=\{\omega \in \{0,1,...,N-1\}^{\bf N}:\forall n\in
{\bf N} , T^n(x)\in B_{\omega (n)}\}.$$

\noindent The set $\varphi _\beta(x)$ is the set of all the possible 
codings of the orbit of $x$ relative to the cover $\beta$.

\begin{definition}
The information content  of $n$ steps of the orbit of $x$ with respect to $\beta$ is defined as $$
{K}(x,T,\beta,n)=\stackunder {\omega \in \varphi _\beta (x)}{min} {AIC_U(\omega ^n)}.$$ $$
\hat{K}(x,T,\beta,n)=\stackunder {\omega \in \varphi _\beta (x)}{min} {AIC_U(\omega ^n|n}).$$

\noindent where $\omega ^n$ is the string containing the first $n$ digits of $\omega $
\end{definition}

We are interested to the asymptotic behavior of this information 
content when $n$ goes to infinity.
We  give a measure of such an  asymptotic behavior by comparing the quantity of information necessary to describe $n$ step of the orbit  with a  function $f$  whose asymptotic behavior is known. For each monotonic function $f(n)$ with $ \stackunder {n\rightarrow \infty }{lim}f(n)=\infty$ we define an indicator
of orbit complexity by comparing the asymptotic behavior of  ${AIC}_U(\omega ^n)$  
or ${AIC}_U(\omega ^n|n)$ with $f$.  
From now on, in the definition of indicators $f$ is always assumed to be monotonic and tends 
 to infinity.

\begin{definition}\label{defin}
The complexity of the orbit of $x\in X$ relative to $f$ and  $\beta$ is defined as:

$$K^f(x,T,\beta)=\stackunder {n\rightarrow \infty}{limsup}
\frac {K(x,T,\beta,n)}{f(n)}
$$

in a similar way we define

$$\hat{K}^f(x,T,\beta)=\stackunder {n\rightarrow \infty}{limsup}
\frac {\hat{K}(x,T,\beta,n)}{f(n)}.
$$

\end{definition}

As it is intuitive, if we refine the cover the information needed to describe 
the orbit increases.
 
\begin{lemma}\label{16}
If $\alpha $ and $\beta $ are open covers of $X$ and $\alpha $ is a refinement of $\beta$ then for all $f$ 
\begin{equation}\label{raffinato}
K^f(x,T,\beta)\leq K^f(x,T,\alpha)
\end{equation}
\begin{equation}\label{+raffinato}\hat{K}^f(x,T,\beta)\leq \hat{K}^f(x,T,\alpha).\end{equation}
\end{lemma}
The proof of Eq. \ref{raffinato} is contained in \cite{ionewnew}
the proof of Eq. \ref{+raffinato} is essentially the same.

As it was said before the definition of $\hat{K}^f$ is particularly suited for
very regular strings. It considers the information contained in the string 
without considering the information contained in its length. Since this quantity is less or equal than $log(n)$ then $|\hat{K}(x,T,\beta,n)-{K}(x,T,\beta,n)|\leq log n +C$ and this implies

\begin{proposition}\label{equu} If $log(n)=o(f(n))$ then $\hat{K}^f(x,T,\beta)=K^f(x,T,\beta)$.
\end{proposition}

 Taking the supremum   over the set of all {\em finite open} covers
$\beta$ of the metric space $X$
 it is possible to get rid of the dependence of our definition  on the choice of the cover $\beta $ and define the complexity of the orbit of
$x$:

\begin{definition}The complexity of $x$ with respect to $f$ is defined as
$$K^f(x,T)=\stackunder {\beta \in \{Open  \ covers\}}{sup}( K^f(x,T,\beta)). $$
$$\hat{K}^f(x,T)=\stackunder {\beta \in \{Open  \ covers\}}{sup}(\hat{ K}^f(x,T,\beta)). $$
\end{definition}

This definition associates to a point belonging to $X$ and a function $f$ a real number which is a measure of the complexity of the orbit of $x$ with respect to the asymptotic behavior of $f$. 
In some sense in the above definition the function $f$ plays a role similar 
to the parameter $d$ in the definition of $d$-dimensional Hausdorff measure.
Each orbit will have a class of functions $f$  such that $K^f(x)$ is finite, 
characterizing the asymptotic behavior of the information that is necessary 
to describe the orbit.

We remark that in the definition above we used two different notions
of ``complexity'' of a string. In principle the construction we made
allows to use any measure of complexity of a string. In \cite{io3}, \cite{CSF} for example
a computable notion of information content based on data compression 
algorithms is considered.

Generalized orbit complexity is invariant under topological conjugation,
as it is stated in the following theorem whose proof follows directly from the definitions:

\begin{theorem}[Invariance] If the dynamical systems $(X,T)$ and $(Y,S)$ are topologically  conjugate,   $\pi:X\rightarrow Y$ is the conjugating homeomorphism, and
$\pi(x)=y$  then ${K}^f(x,T)={K}^f(y,S)$ and  $\hat{K}^f(x,T)=\hat{K}^f(y,S)$.
\end{theorem}

From now on in the notation $\hat{K}^f(x,T)$ we will avoid to explicitly mention the map $T$ when it is clear from the context. 
We now give some example of different behaviors of $\hat{K}^{f}(x)$.

{\em \underline {Periodic orbits.}} If $x$ is a periodic point 
some of the symbolic coding of its orbit is a  periodic string.
An $n$ digit long periodic string can be generated by a program
containing the first period and the length of the string. Since $n$ is given by the definition 
of conditional information content
 $\hat{K}(x,n,V)\leq C$, where $C$ is a constant
not depending on $n$ and $\hat{K}^f(x)=0$ for each $f$.

{\em \underline{Positive entropy.}} 
In the positive entropy case the main result is the following

\begin{theorem}[Brudno's main theorem]\label{brubno}
Let $(X,T)$ be a dynamical system over a compact space. If $\mu$ is an ergodic probability
measure on $(X,T)$, then

$$K^{id}(x,T)=h_{\mu}(T)$$
 
\noindent for $\mu$-almost each $x\in X$ ($id$ is the identity function and $h_{\mu}(T)$ is the KS entropy).
\end{theorem}

By theorem \ref{brubno} and Proposition \ref{equu} it also follows  that if a system is compact, ergodic and has positive Kolmogorov entropy then for almost all points we have  ${\hat K}^{id}(x)=h_\mu$ (and ${\hat K}^{f}(x)=\infty$ if $f=o(id)$).

{\em \underline{ Manneville map.}}
An important example is the piecewise linear Manneville map:
\begin{equation}\label{**}
T_z(x)=\left\{ \begin{array}{cc}
\frac{\xi _{k-2}-\xi _{k-1}}{\xi _{k-1}-\xi _k}(x-\xi _k)+\xi _{k-1} & \xi
_k\leq x<\xi _{k-1} \\ 
\frac{x-a}{1-a} & a\leq x\leq 1
\end{array}
\right. 
\end{equation}

\noindent with $\xi _k=\frac a{(k+1)^{\frac 1{z-1}}}$ ,$\ k\in N,z>2$.
This is a piecewice linear  version of the Manneville map ${\cal T}_z(x)=x+x^z \ (mod\  1)$.

 The Manneville map was introduced in \cite{Manneville}  as an extremely simplified 
model of intermittent behavior in turbulence, then its mathematical 
properties was studied by many authors, And the map was applied as a model of other physical phenomena.
By \cite{ionewnew} which follows from \cite{Gaspard} we have that if $\epsilon$ is small enough then

$$\int _{[0,1]} \hat{K}(x,n,\epsilon) dx\sim n^{\frac z{z-1}}$$

i.e. the Lesbegue average  information that is necessary to describe the orbits of the Manneville map 
for $z>2$ increases as a power law with exponent $\frac 1{z-1}$. 
We remark that the above Manneville map has positive topological entropy 
then it is not weakly chaotic in a topological sense. By the result above we 
can say that the Manneville map in some sense is weakly chaotic with respect 
to the Lesbegue measure.

{\em \underline{ Logistic map at the chaos threshold.}}
Now we calculate the complexity of the orbits of this widely studied  dynamical system. We state a result that, using similar techniques 
slightly improves a result of \cite{menconi} about the complexity of such a
map.

 To understand the dynamic  of the logistic map at the chaos threshold let us use a result
of \cite{collet} (Theorem III.3.5.)

\begin{lemma} \label{teckmann}
The logistic map $f_{\lambda_\infty}$ at the chaos threshold has an
invariant Cantor set $\Omega$. 

\noindent (1) There is a decreasing chain of closed subsets 
$$ J^{(0)} \supset J^{(1)} \supset J^{(2)} \supset \dots, $$ each of
which contains $1/2$, and each of which is mapped onto itself by
$f_{\lambda_\infty}$.

\noindent (2) Each $J^{(i)}$ is a disjoint union of $2^i$ closed intervals. 
$J^{(i+1)}$ is constructed by deleting an open subinterval from the
middle of each of the intervals making up $J^{(i)}$.

\noindent (3) $f_{\lambda_\infty}$ maps each of the intervals making 
up $J^{(i)}$ onto another one; the induced action on the set of
intervals is a cyclic permutation of order $2^i$.

\noindent (4) $\Omega = \cap_i J^{(i)}$. $f_{\lambda_\infty}$ maps 
$\Omega$ onto itself in a one-to-one fashion. Every orbit in $\Omega$ is
dense in $\Omega$. 

\noindent (5) For each $k \in {\bf N}$, $f_{\lambda_\infty}$ has exactly one 
periodic orbit of period $2^k$. This periodic orbit is repelling and
does not belong to $J^{(k+1)}$. Moreover this periodic orbit belongs
to $J^{(k)} \setminus J^{(k+1)}$, and each point of the orbit belongs
to one of the intervals of $J^{(k)}$.

\noindent (6) Every orbit of $f_{\lambda_\infty}$ either lands after 
a finite number of steps exactly on one of the periodic orbits
enumerated in 5, or converges to the Cantor set $\Omega$ in the sense
that, for each $k$, it is eventually contained in $J^{(k)}$. There are
only countably many orbits of the first type.
\end{lemma}

\begin{theorem}
In the dynamical system $([0,1],f_{\lambda _\infty})$, {\em for each} $x\in[0,1]$ 
and each $f$ $\hat{K}^{f}(x)= 0$.
\end{theorem}
{\em Proof.} By the theorem above (point 6) we have that each point $x$ either
is eventually periodic (and the statement follows immediately) either
its orbit converges to the attractor, that is the orbit of 
$x$ is eventually contained in $J_k$ for each $k$.
Now let us consider a cover $V$ and let $\epsilon _V$ its Lesbegue constant.
Let $K_V$ be such that each interval in  $J_k$ has diameter less than $\epsilon _V$. Now if $k>K_V$ each interval of $J_k$ is contained in some set of $V$.
Moreover, by point 3 above we know that the action of the map over $J_k$ is
periodic. This implies that in the set $\varphi _V (x)$ of symbolic orbits 
of $x$ there is an eventually periodic string and then the statement follows easily. $\Box$


{\em \underline{Chaotic maps with zero topological entropy.}}
In \cite{smital} Smital showed an interval map that is continuous, it has 0-topological
entropy and it is chaotic in the sense of Li and Yorke.
The Smital's weakly chaotic maps have non trivial
orbit complexity, that is: there is an uncountable set of points $S$
such that for each $x\in S$ there is $f$ such that $\hat{K}^f(x)>0$.
This is implied by Theorem \ref{mitico} below.

In order to prove it we give the definition of weak scattering set, this is a notion
that is weaker than the notion of scattering set used in the Li-Yorke
definition of chaotic map (see \cite{smi2} e.g.).

\begin{definition}
A set $S$ is called a weak scattering set if there is a $\delta $ such that for all $x,y\in S, x\neq y$ 
implies $\stackunder {n\rightarrow \infty}{limsup} d(T^n(x),T^n(y))>\delta$.
\end{definition}
 
A point is said to have nontrivial orbit complexity if $\exists f$ with $\hat{K}^f(x)>0$.
The following is an easy criterion to find systems with nontrivial
orbit complexity.

\begin{theorem}\label{mitico}
If $(X,T)$ has an uncountable weak scattering set $S$, then
there is a set $S'$ with $\#(S')<\#(S)$ such that for each $x\in S-S'$ 
 $x$ has nontrivial orbit complexity.
\end{theorem}

{\em Proof.} The proof is based on a cardinality argument.
First we proof that $\exists x,f$ such that $K^f(x)>0$.
Conversely let us suppose that there are not such points.
This implies that given any cover $V=\{ B_1,...,B_v\}$ for each $x$ there is a finite 
set of programs $P_x=\{{p^1}_x,...,{p^k}_x\}$ such that $\forall n$  the simplest symbolic orbit
 $\omega ^n$ of $x$ with respect to $V$  (that is such that $AIC(\omega ^n|n)=\stackunder{\omega \in \varphi _V (x)}{min}AIC(\omega ^n|n)$) is such that $\omega ^n ={\cal U}({p^i}_x, n)$ for some $ {p^i}_x\in P_x$. Now the set $P=\{P_x|x\in X\}$ is countable
because is contained in the finite parts of a countable set (the set of all
programs). This implies that there is a set $Z\subset S$ with $\#(Z)=\#(S)$ such that for each $x\in Z$ the set $S_x=\{ y\in S, P_y=P_x\}   $
is uncountable.  
 
Now  the set of possible infinite symbolic orbits
associated to the set of programs $P_x$ is finite.
Let us consider  the set $W_{P_x}=\{ \omega\in \{1,...,v\}^{\bf N} \ s.t. 
\forall n \exists i\leq k , \omega ^n = {\cal U}({p_x}^i,n)\}$
(that is the set of symbolic orbits that can be generated from the set of programs $P_x$) this set is finite and has $\# (W_{P_x}) \leq k$.
This leads to a contradiction. By the definition of $W_{P_x} $ each point of $S_x $ must be such that there is an $\omega \in W_{P_x}\cap \varphi _V (x)
$ (the set  $W_{P_x} $ is the set of the possible orbits related to $P_x$).
But now, since $S_x\subset S$ is weakly scattering if $diam(V)<\delta $
 we have that if $x\neq y$ then $ \varphi _V (x)\cap \varphi _V (y)=\emptyset$
because at some time the orbits of the two points will be contained in  sets of $V$ having  of course empty intersection because the distance of the two points is such that they cannot be in the same set of $V$.
Since $W_{P_x}$ is finite then is not possible that $\forall x \ W_{P_x}\cap \varphi _V (x)\neq \emptyset$. This ends the first part of the proof.
Now let us consider the scattering set $S_1$ with $S_1=S-x$. $S_1$ does not contain $x$  and still verifies the assumptions
of the theorem, then by the first part of this proof there is $y\neq x$ such that $y$ has nontrivial orbit complexity. In this way by induction we prove the full theorem.$\Box$

\begin{remark}\label{misi} By a result of 
Misiurewikscz (\cite{misiu}) each zero entropy, continuous map on the interval is metrically  hisomorphic to the adding machine for each non atomic invariant measure.

Then even the Smital's map are  hisomorphic to the adding machine from 
the measure preserving point of view. 
The orbit complexity of such maps is different from the complexity
of other zero entropy continous maps of the interval (logistic map  e.g.). 

This implies
that orbit complexity in this case is more sensitive than any invatiant
 constructed in the measure-preserving framework.
\end{remark}

\section{Generalized topological  entropy}

We give the definition of an indicator of the global topological complexity of the system. 
Similarly to the classical definition of
topologial entropy.
This indicator will measure the asymptotical
behavior of the number of substantially different orbits that appears
in the dynamics.

Let $X$ be a compact metric space and $T\in C^0(X,X)$.
If $x,y\in X$ let us say that $x,y$ are $(n,\epsilon)$ separated if  $d(T^k(x),T^k(y))>\epsilon$
for some $k\in\{0,...,n-1\}$. If $d(T^k(x),T^k(y))\leq \epsilon$ for each $k\in\{0,...,n-1\}$ then $x,y$ are said to be $(k,\epsilon)$ near. 
A set $E\subset X$ is called $(n,\epsilon)$ separated if $x,y\in E,x\neq y$  
 implies that $x$ and $y$ are $(n,\epsilon) $ separated.

Moreover a set $E\subset X$ is called an  $(n,\epsilon)$ net if $\forall x \in X $  
  $\exists y\in E \ s.t. $ $x$ and $y$ are $(n,\epsilon)$ near.

We remark that a maximal $(\epsilon, n)$ separated set is a $(2\epsilon,n) $ net.

Let us consider as in the classical definition of topological entropy

$$ s(n,\epsilon)=max\{card( E):E\subset X \ is \ (n,\epsilon)-separated\}.$$

We choose a monotone function $f$  such that $\stackunder {n\rightarrow \infty}{lim} f(n)=\infty $ as before and define

$$h^f(T,\epsilon)=\stackunder{n\rightarrow \infty}{limsup} \frac{\log (s(n,\epsilon))}{f(n)}$$
$h^f(T,\epsilon)$ is monotone in $\epsilon$ so we can define
$$h^f(T)=\stackunder{\epsilon \rightarrow 0} {lim}h^f(T,\epsilon).$$

Let us consider 
$$r(n,\epsilon )=min\{ card( E):E\subset X ,E \ is \ a \ (n,\epsilon ) \ net\}.
$$

Since $r(n,\epsilon) \leq s(n,\epsilon)$ and (see \cite{petersen} p. 268) $s(n,\epsilon)\leq r(n,\frac \epsilon 2)$
in the definition of generalized topological entropy $h^f$ we can also consider instead of  $s(n,\epsilon)$ the number
$r(n,\epsilon )$ and obtain an equivalent definition.

We also give a third equivalent possible definition, using open covers.
Let $U $ and $V$ two open covers of $X$. We denote by $U\vee V$ the least common refinement  
(or join) of $U$ and $V$.
If $U$ is an open cover, let us denote by $N(U)$ the minimum cardinality
of the subcovers of $U$, and let $H(U)=log(N(U))$.
Let us consider $${\tt h}^f(T,U)=\stackunder {n\rightarrow \infty}{limsup}\frac {H(U\vee T^{-1}(U)\vee...\vee T^{-n}(U))}{f(n)}$$
and $${\tt h}^f(T)=\stackunder {finite\ open\ covers}{sup}{\tt h}^f(T,U).$$
We remark that   (see \cite{petersen} page 268) 
if $\epsilon > diam(U)$ then

$${N(U\vee T^{-1}(U)\vee...\vee T^{-n}(U))}\geq s(n,\epsilon)$$
and if $\epsilon $ is the Lesbegue number of  $U$ then

$${N(U\vee T^{-1}(U)\vee...\vee T^{-n}(U))}\leq r(n,\epsilon).$$

By this it follows that

\begin{proposition}\label{10} For each $f$ ${\tt h}^f(T)= h^f(T)$ and    
the definitions are equivalent.
\end{proposition}

By this is also follows that $h^f$ is invariant under 
isomorphisms of dynamical systems.

\begin{proposition} If $(X,T)$, $(Y,T')$ are topological dynamical systems, $\psi:X\rightarrow Y $ is an
homeomorphism such that 

\begin{equation}\label{equisopra}
\begin{array}{rcccl}
\ & \ & \psi & \ & \ \\
\ & X & \rightarrow & Y \ \\
T & \downarrow & \ & \downarrow & T'\\
\ & X & \rightarrow & Y \ \\
\ & \ & \psi & \ & \ \\
\end{array}
\end{equation}
commutes, then for each $f$ $h^f((X,T))=h^f((Y,T'))$.
\end{proposition}

{\em Proof. }
Since $\psi$ is an homeomorphism then it sends an open cover of $X$ to 
an open cover of $Y$ and it is easy to see that since the diagram 
commutes ${\tt h}^f((X,T))={\tt h}^f((Y,T'))$. By Proposition \ref{10} we also have ${ h}^f((X,T))={h}^f((Y,T'))$. $\Box$ 

Now we state a result which is useful to characterize the systems where
the generalized topological entropy is null for each $f$.

\begin{definition} A system $(X,T)$ is 
said to be {\it
equicontinuous} if for any $\epsilon >0$ there is $\eta>0$ such that 
if $x,\ y \in X$ with $d(x,y) < \eta$ then for any $n \in {\bf N}$ one has 
$d(T^n(x), T^n(y)) < \epsilon$.
\end{definition}
 If $(X,T)$ is not 
equicontinuous there are $\epsilon
>0$ and a point $x \in X$ such that for any $\eta > 0$ one can find $y \in X$
with $d(x,y)<\eta$ and $n \in {\bf N}$ such that $d(T^n (x),T^n (y)) > 
\epsilon$.

 A point $x \in X$ is called an {\it equicontinuity point} 
if for any $\epsilon >0$ there is $\eta>0$ such that if $y \in X$ 
with $d(x,y) < \eta$ then for any $n \in {\bf N}$ one has $d(T^n(x), T^n(y)) < 
\epsilon$; obviously a system is equicontinuous if all its points are 
equicontinuity points. 



Let us cite the following fundamental result from \cite{blanc} about the topological complexity  of a dynamical system.

\begin{theorem}\label{bla} Let $(X,T)$ be a dynamical system, $X$ is compact, $T$ is continuous The two following statement 
are equivalent:

{(1)} $(X,T)$ is equicontinuous.

{(2)} For any finite open cover ${ U}$ of $X$, $N( U\vee T^{-1}(U)\vee ...\vee T^{-n}(U))$ is bounded.
\end{theorem}

By this it follows that

\begin{proposition}
 $(X,T)$  is  equicontinuous if and only if for each $f$ it holds $ h^f((X,T))=0$.  
\end{proposition}
 
The next result shows that for the logistic map at the chaos threshold   $\exists$ $f$ such that $\hat{K}^f(x)<h^f(T)$ for each $x\in [0,1]$. 
This is an example of a system where local and global complexity are quite different. 
In the following (Theorem \ref{infinit}) we will see  that in general  $\hat{K}^f(x)\leq h^f(T)$ for a wide class  
of dynamical systems.

\begin{theorem}The  logistic map $f_{\lambda _\infty}$
 is not equicontinuous and then there is some $f$ such that $h^f(f_{\lambda _\infty})>0$.
\end{theorem}
{\em Proof.} By Lemma \ref{teckmann} we have that either an orbit
is eventually periodic or it converges to  the attractor $\Omega$ and the orbit is dense on the attractor, then let $p$ be periodic and  $r=d(orb(p),\Omega)$. In each neighborhood of the periodic point there is a point $x$ that converges to the attractor and then
 $limsup( d(T^n(x),T^n(p)))>r$. The point $p$ cannot be an  equicontinuity point. The statement follows then from Theorem \ref{bla} $\Box$

\section{Computable  Structures, Constructivity }\label{comp}

In this section we give a rigorous notion of constructive map and
the results about constructive mathematics that are necessary in the following. Constructive functions can be considered in some sense  as  algorithms acting over metric spaces.
All function that can be concretely defined and effectively calculated
are constructive.
Algorithms works with strings, if strings are interpreted as points
of a metric space we have the possibility to relate the world of continuum 
mathematics with the world of algorithms.
This is what is currently done when expressing a point of a metric space
by a symbolic notation. For example $\frac \pi 2$ is a symbolic string
that represents a point of the metric space ${\bf R}$ and  allows to calculate a symbolic representation for the value of
$sin(\frac \pi 2) $ by some algorithm (because $sin $ is a constructive function).
   An interpretation function is a way to interpret  a string as a point of the metric space.
 An interpretation is said to be computable  if the distance between ideal points is computable with arbitrary precision:

\begin{definition}[Computable interpretation]
A computable  interpretation function on  $(X,d)$ is a function $I:\Sigma \rightarrow
X$ such that $I(\Sigma)$ is dense in $X$ and  there exists a total recursive
function $D:\Sigma \times \Sigma \times
 {\bf N} \rightarrow {\bf Q}$ such that $\forall s_1,s_2 \in \Sigma ,n\in
{\bf N}$:

$$|d( I(s_1),I(s_2))-D(s_1,s_2,n)|\leq \frac{1}{2^n}.$$

\end{definition}

 A point $x\in X$ is said to be {\em ideal} if it is the image of some string: $x=I(s), s\in \Sigma$.  

Two interpretations are said to be equivalent if the distance from an ideal
point from the first and a point from the second is computable up to
arbitrary precision.

\begin{definition}[Equivalence of interpretations]\label{equi7}   
Let $I_1$ and $I_2$ be two computable interpretations in $(X,d)$;
we say that $I_1$ and $I_2$ are equivalent if there exists a total recursive function
$D^*:\Sigma \times \Sigma \times {\bf N} \rightarrow {\bf Q}$,
such that $\forall s_1,s_2 \in \Sigma ,n{\in {\bf N}}$:

$$|d( I_1(s_1),I_2(s_2))-D^*(s_1,s_2,n)|\leq \frac{1}{2^n}.$$

\end{definition}

For example,  finite binary  strings $s\in \Sigma $ can be interpreted as rational numbers by interpreting the string as the binary expansion of the number. Another interpretation can be given by 
interpreting a string as an encoding of a couple of integers whose ratio 
gives the rational number. If the  encoding is recursive, the two interpretation are equivalent.

\begin{proposition}
The relation  defined by definition \ref{equi7} is an equivalence relation.
\end{proposition}
For the proof of the above proposition see \cite{io1}.

\begin{definition}[Computable structure]
A computable structure ${\cal I}$ on $X$ is an equivalence class  of computable  interpretations  in $X$.
\end{definition}

%
%
\begin{remark}\label{tondo} 
We remark as a property of the computable structures that if $B_r(I(s))$ is 
an open ball with center in an ideal  point $I(s)$ and rational radius $r$ and $I(t)$ is another point   then there is an algorithm that 
verifies if  $I(t)\in B_r(I(s)) $. If  $I(t)\in B_r(I(s)) $ then the algorithm
outputs ``yes'', if $I(t)\notin B_r(I(s)) $ the algorithm  outputs ``no'' or 
does not stop. The algorithm calculates $ D(s,t,n)$ for each $n$ until it
finds that $D(s,t,n)+2^{-n}< r$ or  $D(s,t,n)-2^{-n}> r$, in the first case
it outputs ``yes'' and in the second it outputs ``no''. If $d(I(s),I(t))\neq r$ the algorithm will stop and output an answer.
\end{remark} 

We give a definition of {\em morphism} of metric spaces with computable 
structures. A morphism is heuristically  a computable
function between  computable metric spaces.
 The definition states that if $\Psi$ is a morphism the image  
of an ideal point can be calculated up to arbitrary precision by an algorithm.

\begin{definition}[Morphism between computable structures]
If $(X,d,{\cal I})$ and $(Y,d',{\cal J}) $ are spaces with  computable
structures; a function $\Psi :X\rightarrow Y$ is said to be  a morphism of
computable structures if $\Psi $ is  continuous and for each pair
$I\in {\cal I},J\in {\cal J}$ 
there exists a total recursive function
$D^*:\Sigma \times \Sigma \times {\bf N} \rightarrow {\bf Q}$,
such that $\forall s_1,s_2 \in \Sigma ,n{\in {\bf N}}$:

$$|d'(\Psi( I(s_1)),J(s_2))-D^*(s_1,s_2,n)|\leq \frac{1}{2^n}.$$
\end{definition}

\noindent We remark that $\Psi $ is not required to have dense image
and then $\Psi(I(*))$ is not necessarily an interpretation function 
equivalent to $J$.

\begin{remark}\label{remark10} 
As an example of the  properties of the morphisms, we remark that if a map $\Psi:X\rightarrow Y $  is a morphism  then given a point $x\in I(\Sigma) \subset X$ it is possible to find by an algorithm a point $y\in J(\Sigma) \subset Y$ as near as we want to $\Psi (x) $. 
\end{remark}
The procedure is simple: if $x=I(s)$ and we want to find a point $y=J(z_0)$ such that $d'(\Psi( I(s)),y)\leq 2^{-m} $ then we calculate $D^*(s,z,m+2)$ for each $z\in \Sigma$ until we find $z_0$ such that $D^*(s,z_0,m+2)<2^{-m-1}$. Clearly $y=J(z_0)$ is such that $d'(\Psi(x),y)\leq 2^{-m} $. The existence of such a $z_0$ is assured by the density of $J$ in $Y$.
We also remark that by a similar procedure, given a point $I(s_0)$ and $\epsilon \in {\bf Q}$ it is possible to find a point $I(s_1)$ such that $d(I(s_0),I(s_1))\geq \epsilon $.

A constructive map is a morphism for which the continuity relation between $\epsilon$ and $\delta$  is given by a
 recursive function.

\begin{definition}[Uniformly constructive functions]\label{definition11}
A function

 $\Psi:X\rightarrow Y $ between  spaces with computable structure
$ (X,d,{\cal I}) $, $(Y,d',{\cal J})$ is said to be uniformly constructive
if  $\Psi$ is a morphism between the  computable structures and it is 
effectively uniformly continuous, i.e. there is a total recursive function
$f:{\bf N}\rightarrow {\bf N} $ such that for all $x,y \in X$  
$d(x,y)<2^{-f(n)}$ implies $d'(\Psi (x),\Psi (y))<2^{-n}$.
\end{definition}

If a map between spaces with a computable structure is uniformly constructive then there
is an algorithm to follow the orbit  each ideal point $x=I(s_0)$.

\begin{lemma}\label{fascino}
If $T: (X,{\cal I})\rightarrow (X,{\cal I})$  is uniformly constructive, $ I\in{\cal I}$ then  there is an algorithm (a total recursive function)
$A:\Sigma \times {\bf N}\times {\bf N}\rightarrow \Sigma$ such that $\forall k,m\in {\bf N},s_0\in \Sigma$ $d(T^k(I(s_0)),I(A(s_0,k,m)))< 2^{-m} $. 
\end{lemma}
{\em Proof.}
Since $T$ is  effectively uniformly continuous
we define the function $g_k(m)$ inductively as 
$g_1(m)=f(m )+1$, $g_i(m)=f(g_{i-1}(m)+1)$ where $f$ is the function of effective uniform
continuity of $T$ (definition \ref{definition11}).
 If  $d(x,y)<2^{-g_k(m)} $ then $d(T^i(y),T^i(x))<2^{-m}$ for $i\in \{1,...,k\}$. 
Let us choose $I \in {\cal I}$. 
We recall that the assumption that $T$ is a morphism implies that
there is a recursive function $D^* (s_1,s_2,n)$ such that 

$$|D^* (s_1,s_2,n)-d(I(s_1),T(I(s_2)))| <2^{-n}.$$

\noindent Let us suppose that $x=I(s_0)$.
Now let us describe the algorithm $A$:
 using the function  $D^*$ and the function $f$,  $A$ calculates $g_k(m)$ and finds a string $s_1$ such that 
$d(I(s_1),T(I(s_0))) <2^{-g_k(m)}$ as described in remark \ref{remark10}. This is the first step of the
algorithm. Now $d(T(I(s_1)),T^2(x))\leq 2^{-(g_{k-1}+1)}$.
We can use $D^*$ to find a string $s_2$ such that $d(I(s_2),T(I(s_1)))<2^{-(g_{k-1}+1)}$. By this $d(I(s_2),T^2(x)) \leq 2^{-g_{k-1}}$. 
 This implies that $d(T(I(s_2)),T^3(x))\leq 2^{-(g_{k-2} +1)}$, then 
we find $s_3$ such that $d(I(s_3),T(s_2))\leq 2^{-(g_{k-2} +1)}$ and so on for
$k$ steps. At the end we find a string $s_k$ such that $d(I(s_k),T^k(x)) \leq 2^{-m}$. $\Box $

Let us  describe a certain class of {\em nice} balls covers that will be used in the following. 
\begin{definition} If ${\alpha } =\{B_1(y_1,r_1),...,B_n(y_n,r_n)\}$ is a
ball cover of the metric space $X$  whose elements are balls with centers $y_i$
and radii $r_i$. We say that $\alpha$ is a {\em nice}
cover if $X\subset \stackunder {i}{\cup}B_i(y_i,\frac {r_i}{2} )$.
\end{definition}
In other words  ${\alpha }$ is a {\em nice } cover if  dividing the radius of
the balls by $2$ we  have  again  a cover.

\begin{remark}
We remark that since the space is compact each open cover has a refinement
which is a nice cover.
\end{remark}

\begin{remark}If we have a nice cover ${\alpha } =\{B_1(y_1,r),...,B_n(y_n,r)\}$ made of balls with ideal centers $y_i=I(s_i)$ and rational radius and we have $x\in X$ and an ideal point $y=I(s)$ such that $d(x,y)<\frac r2$ then it is possible to find a ball of $\alpha$ that contains $x$. Indeed we find 
by the properties of computable interpretations an $y_i$ such that $d(y,y_i)<\frac r2 $ (see Remark \ref{tondo}). This is possible because the cover is nice. The ball in the cover with center $y_i$ will contain $x$.
\end{remark}

By the above remark if we have an algorithm $A()$ to follow the orbit 
of  ideal points as in Lemma \ref{fascino} and a nice cover ${\alpha } =\{B_1(y_1,r),...,B_n(y_n,r)\}$ made of balls with ideal centers and rational radius 
 for each ideal point $x$ and $n$ it is possible to find by an algorithm a set $B$ of $\alpha $ such that $T^n(x)\in  B$ and then
we can construct a symbolic orbit of $x$ relative to $\alpha$.
We remark that if $n$ is given (by the definition of $AIC(s|n)$) the lenght of a code implementing such 
an algorithm does not depend on $n$.
Then we have
\begin{lemma}\label{magia}
Let $(X,{\cal I})$ is a metric space with a computable structure and
\begin{itemize}
\item $\{B(x_i,r)\}_{1\leq i\leq m} $ is a nice cover made of balls
with ideal centers and rational radius $r$,

\item $x\in X$ is an ideal point

\item $P_{\frac r2}:{\bf N} \rightarrow \Sigma$ is such that $\forall k\in \{1,...,n\}$  
$$d(T^k(x),I(P_{\frac r2}(k)))<\frac r2 $$
is an algorithm to follow the orbit of $x$ with  accuracy $\frac r2 $ similar as above in Lemma \ref{fascino}. 
\end{itemize}
then there is an algorithm $P':{\bf N} \rightarrow \{1,...,m\}$ such that $\forall i\leq k, P'(i)=j $ if $T^i(x)\in B(x_j,r)$.
Moreover the length of a code implementing $P$ on an universal 
Turing machine is equal to the code for $P'$ up to a constant that 
does not depend on $n$. 
\end{lemma}

\section{Gen. Top. Entr. and Orbit complexity}

\begin{theorem}\label{infinit}
If $X$ is compact and $(X,T)$ is constructive (for some computable
structure over $X$) then 
for each $f$ and each $x\in X$ 
$$ \hat{K}^f(x)\leq h^f(T).$$
\end{theorem}
{\em Proof.}
Since the system is constructive
by the use of the algorithm $A()$ (Lemma \ref{fascino}) following the orbit of an ideal point
at any given accuracy and $D()$ approximating the distance $d$ at any 
given accuracy we have the following.
If $x,y$ are ideal, $I(z)=x, I(w)=y $  and $k\in {\bf N},\epsilon \in {\bf Q}$ 
then there is a procedure 
 $P(z,w,\epsilon,k)$ such that if there is a $i\leq k$ such that $d(T^i(x),T^i(y))>\epsilon$  then the procedure stops with output ``YES''.

We also remark that if $x=I(z)$ is ideal $I(z)=x$ and there is $y\in X$ such that 
$x,y$ are $(n,\epsilon) $ separated then we can find by a procedure
 $P'(z,n,\epsilon)$ an ideal $y'$ such that $x,y'$ are $(n,\epsilon) $ separated.
The procedure $P'$ calculates $P(z,v,n,\epsilon)$ for all strings $v$
in a ``parallel'' way until it finds a $v$ stopping the procedure
with output YES. Such a $v$ must exist by the density of $I$ and 
the continuity of the map $T$.


We will prove that 
if $\beta=\{B(x_1,2\epsilon),...,B(x_n,2\epsilon)\}$
 is a nice cover with ideal centers and rational radius and $\epsilon $ is small enough, then
 for each $k\in {\bf N}$ and $\forall x \in X$ there is  a program $p_{k,\epsilon}$ such that ${\cal U}(p_{k,\epsilon},k)$ is a symbolic coding of the first $k$ steps 
of the orbit of $x$ with respect to $\beta $ and 
$|p_{k,\epsilon}|\leq  log(s(k,\frac \epsilon 2))+C $.

The idea is that by the constructivity (as it is said above) we can construct $(n,\epsilon) $ 
separated sets and use them to select the points that give rise to 
separate orbits, moreover the number of these points is bounded by 
 the topological entropy.
Now let us see a more precise description of such a  procedure.

The program $p_{k ,\epsilon}$ we want to describe now  will contain a number $n_x$ with  $n_x\leq s(k,\frac \epsilon 2)$ and a procedure to construct a symbolic orbit of $x$.
The procedure runs as follows 

 \ \ \ {\tt First let us consider an empty list
of strings: list=$\emptyset $

\

 \ \ \  For each $i\in {\bf N},i\leq n_x$ do :$\{$ by the above procedure $P'()$ find an ideal point
 that is $(k,\frac{\epsilon}2 ) $ separated from the points in the list
 and add the corresponding string to the list $\}$

\

 \  \ \ Follow by the algorithm $A()$ (Lemma \ref{fascino}) with accuracy $\epsilon $
 the last point found in the list. 

\ \ \ By the procedure stated in Lemma \ref{magia}  produce a symbolic string
 associated to the cover $\beta$.

\ \ \  
}

We remark that in the above procedure the number $k$ is given 
by the definition of $AIC(s|k)$.

%
If $n_x$ is big enough
the previous procedure construct a maximal $(k,\frac \epsilon 2)$ separated
set.
Now since a maximal $(k,\frac {\epsilon} 2)$ separated set is also a $(k,\epsilon)$
net, then each $x\in X$  is $ \epsilon$ near to some point 
found by the procedure above.
The procedure then  follows by the algorithm $A()$ with accuracy $\epsilon $
 the orbit of such a point and then produce the symbolic list associated with
the given cover. Since   $n_x\leq s(k,\frac \epsilon 2)$ then $|p_k|\leq log( s(k,\frac \epsilon 2)) +C$ and the statement follows.
$\Box $  

The assumption of constructivity in the above theorem is essential.
 There are examples of dynamical systems that are not constructive
 for any computable structure and violates the above inequality for 
each point.

\begin{example}
Let us consider $(X,T)$ with $X=S^1$ and $T(x)=x+r \ (mod \ 1)$ 
where $r$ is a non constructive number (See \cite{ionewnew} for the definition of such numbers, see below for an example). 

In this example, since $T$ is an isometry we have $h^f(T)=0$ for all $f$.
But if $r$ is not constructive there is some $f$ such that $\hat{K}^f(x)>0$ for 
all $x$. 
\end{example}

The idea is that the knowledge of many steps of  a symbolic orbit for $x$ with respect 
to the cover $\beta$ implies the knowledge of many digits of $r$.

In the following proposition we prove this fact when $r=0.r_1r_2...$
(binary expansion of $r$) is such that $\stackunder {n\rightarrow \infty }{lim } \frac {AIC(r_1...r_n)}n=1$ we call {\em random} such a real (this condition is satisfied for Lesbegue
almost all reals). The proposition in the other cases is a straightforward
generalization.

\begin{proposition}In the above example if $r$ is random we have that for each $x\in X$  $\hat{K}^{log}(x)>0$ and ${h}^f(T)=0$ for each $f$.
\end{proposition}
{\em Proof.} 
Let $x\in S^1$ and $\beta=\{B(x_1,r),...,B(x_n,r) \}$
a cover of $S^1$ with rational centers and radius. Without loss of generality we can suppose $x=0$. 

 Let $p_k$ be a program generating a sequence $\omega_1,...,\omega_{k -1} $
 such that $T^i(0)\in B(x_{\omega _i})$.

By the use of $p_k$ it is possible to find a rational $q$ such that in the dynamical 
system $T':S^1\rightarrow S^1 $ defined by $T'(x)=x+q \ (mod \ 1 )$
the orbit of $0$ has the same associated $k$ steps symbolic orbit.

Since $\forall i\leq k, d(T^i(0),{T'}^i(0))\leq 2\epsilon $ then $|r-q|\leq\frac {2 \epsilon}  k$ and then the knowledge of $k$ steps of the orbit of $0$ with
accuracy $\epsilon$ implies the knowledge of $r$ up to accuracy $\frac {2 \epsilon} 
k$. This implies the knowledge of $log(k) +C$ binary digits of $r$. If $r$ was 
a random real this implies the statement.
$\Box $

In \cite{io2} an example was given of a non constructive system over a non
compact space having large orbit complexity while the map defining it is
equicontinuous. This example showed in the non compact case that constructivity is essential to relate
the complexity  of the behavior of a system and chaos.

The last proposition implies that even in the compact case this is true.
 Constructivity comes out naturally when considering definitions of complexity
 which are based on the algorithmic information content.

\end{document}